\input harvmac.tex

%% MACROS

%% Theorems, Definitions, etc.

\def\thm#1{\bigskip\noindent{\bf Theorem #1} }

\def\rmk#1{\bigskip\noindent{\bf Remarks} }
%%

% Something to deal with sub-sub-sections

\def\unlockat{\catcode`\@=11}
\def\lockat{\catcode`\@=12}

\unlockat
% Something to deal with sub-sub-sections

\def\newsec#1{\global\advance\secno by1\message{(\the\secno. #1)}
\global\subsecno=0\global\subsubsecno=0\eqnres@t\noindent
{\bf\the\secno. #1}
\writetoca{{\secsym} {#1}}\par\nobreak\medskip\nobreak}
\global\newcount\subsecno \global\subsecno=0
\def\subsec#1{\global\advance\subsecno
by1\message{(\secsym\the\subsecno. #1)}
\ifnum\lastpenalty>9000\else\bigbreak\fi\global\subsubsecno=0
\noindent{\it\secsym\the\subsecno. #1}
\writetoca{\string\quad {\secsym\the\subsecno.} {#1}}
\par\nobreak\medskip\nobreak}
\global\newcount\subsubsecno \global\subsubsecno=0
\def\subsubsec#1{\global\advance\subsubsecno by1
\message{(\secsym\the\subsecno.\the\subsubsecno. #1)}
\ifnum\lastpenalty>9000\else\bigbreak\fi
\noindent\quad{\secsym\the\subsecno.\the\subsubsecno.}{#1}
\writetoca{\string\qquad{\secsym\the\subsecno.\the\subsubsecno.}{#1}}
\par\nobreak\medskip\nobreak}

\def\subsubseclab#1{\DefWarn#1\xdef
#1{\noexpand\hyperref{}{subsubsection}%
{\secsym\the\subsecno.\the\subsubsecno}%
{\secsym\the\subsecno.\the\subsubsecno}}%
\writedef{#1\leftbracket#1}\wrlabeL{#1=#1}}% Macros for boxes
\lockat

\input epsf
% Macros plagiarized from P.G.

\def\figin{\epsfcheck\figin}\def\figins{\epsfcheck\figins}
\def\epsfcheck{\ifx\epsfbox\UnDeFiNeD
\message{(NO epsf.tex, FIGURES WILL BE IGNORED)}
\gdef\figin##1{\vskip2in}\gdef\figins##1{\hskip.5in}% blank space instead
\else\message{(FIGURES WILL BE INCLUDED)}%
\gdef\figin##1{##1}\gdef\figins##1{##1}\fi}
\def\DefWarn#1{}
\def\figinsert{\goodbreak\midinsert}
\def\ifig#1#2#3{\DefWarn#1\xdef#1{fig.~\the\figno}
\writedef{#1\leftbracket fig.\noexpand~\the\figno}%
\figinsert\figin{\centerline{#3}}\medskip\centerline{\vbox{\baselineskip12pt
\advance\hsize by -1truein\noindent\footnotefont{\bf Fig.~\the\figno:} #2}}
\bigskip\endinsert\global\advance\figno by1}

\noblackbox

\def\IL{\relax{\rm I\kern-.18em L}}
\def\IH{\relax{\rm I\kern-.18em H}}
\def\IR{\relax{\rm I\kern-.18em R}}
\def\IC{\relax\hbox{$\inbar\kern-.3em{\rm C}$}}
\def\IT{\relax\hbox{$\inbar\kern-.3em{\rm T}$}}
\def\IZ{\relax\ifmmode\mathchoice
{\hbox{\cmss Z\kern-.4em Z}}{\hbox{\cmss Z\kern-.4em Z}}
{\lower.9pt\hbox{\cmsss Z\kern-.4em Z}} {\lower1.2pt\hbox{\cmsss
Z\kern-.4em Z}}\else{\cmss Z\kern-.4em Z}\fi}

\def\CN {{\cal N}}

\def\CB {{\cal B}}

%% MORE MACROS

\def\CN {{\cal N}}

\font\manual=manfnt \def\dbend{\lower3.5pt\hbox{\manual\char127}}

\def\IZ{\relax\ifmmode\mathchoice
{\hbox{\cmss Z\kern-.4em Z}}{\hbox{\cmss Z\kern-.4em Z}}
{\lower.9pt\hbox{\cmsss Z\kern-.4em Z}} {\lower1.2pt\hbox{\cmsss
Z\kern-.4em Z}}\else{\cmss Z\kern-.4em Z}\fi}
\def\half {{1\over 2}}

\def\p{\partial}

\def\CN {{\cal N}}

% more macros, alphabetically

\def\IZ{\relax\ifmmode\mathchoice
{\hbox{\cmss Z\kern-.4em Z}}{\hbox{\cmss Z\kern-.4em Z}}
{\lower.9pt\hbox{\cmsss Z\kern-.4em Z}} {\lower1.2pt\hbox{\cmsss
Z\kern-.4em Z}}\else{\cmss Z\kern-.4em Z}\fi}
\def\IB{\relax{\rm I\kern-.18em B}}
\def\IC{{\relax\hbox{$\inbar\kern-.3em{\rm C}$}}}
\def\ID{\relax{\rm I\kern-.18em D}}
\def\IE{\relax{\rm I\kern-.18em E}}
\def\IF{\relax{\rm I\kern-.18em F}}
\def\IG{\relax\hbox{$\inbar\kern-.3em{\rm G}$}}
\def\IGa{\relax\hbox{${\rm I}\kern-.18em\Gamma$}}
\def\IH{\relax{\rm I\kern-.18em H}}
\def\II{\relax{\rm I\kern-.18em I}}
\def\IK{\relax{\rm I\kern-.18em K}}
\def\IP{\relax{\rm I\kern-.18em P}}
\def\IQ{\relax\hbox{$\inbar\kern-.3em{\rm Q}$}}

\def\inbar{\,\vrule height1.5ex width.4pt depth0pt}

\def\thm#1{\bigskip\noindent{\bf Theorem #1}}

\def\rmk#1{\bigskip\noindent{\bf Remarks #1}}

\def\p{\partial}

\font\cmss=cmss10 \font\cmsss=cmss10 at 7pt
\def\IR{\relax{\rm I\kern-.18em R}}

% Macros for boxes

\def\boxit#1{\vbox{\hrule\hbox{\vrule\kern8pt
\vbox{\hbox{\kern8pt}\hbox{\vbox{#1}}\hbox{\kern8pt}}
\kern8pt\vrule}\hrule}}
\def\mathboxit#1{\vbox{\hrule\hbox{\vrule\kern8pt\vbox{\kern8pt
\hbox{$\displaystyle #1$}\kern8pt}\kern8pt\vrule}\hrule}}

%% ANOTHER SET OF MACROS

\def\inbar{\,\vrule height1.5ex width.4pt depth0pt}

\def\p{\partial}

\font\cmss=cmss10 \font\cmsss=cmss10 at 7pt
\def\IR{\relax{\rm I\kern-.18em R}}

%% new macros

%% END MACROS
%%

\lref\bateman{A. Erdelyi et. al. , {\it Higher
Transcendental Functions, vol. I, Bateman manuscript project}
(1953) McGraw-Hill}

\lref\borchaut{R. Borcherds,
``Automorphic forms with singularities on Grassmannians,''
alg-geom/9609022}

\lref\donint{S.K. Donaldson, ``Connections,
Cohomology and the intersection forms of 4-manifolds,'' J. Diff.
Geom. {\bf 24 } (1986)275.}

\lref\finstern{R. Fintushel and R.J. Stern,
``The blowup formula for Donaldson invariants,'' alg-geom/9405002;
Ann. Math. {\bf 143} (1996) 529.}

\lref\mw{G. Moore and E. Witten, ``Integration over
the $u$-plane in Donaldson theory," hep-th/9709193; Adv. Theor.
Math. Phys. {\bf 1} (1998) 298.  }

\lref\swi{N. Seiberg and E. Witten,
``Electric-magnetic duality, monopole condensation, and confinement
in ${\cal N}=2$ supersymmetric Yang-Mills Theory,'' hep-th/9407087;
Nucl. Phys. {\bf B426} (1994) 19;  ``Monopoles, duality and chiral symmetry
breaking in ${\cal N}=2$
supersymmetric QCD,'' hep-th/9408099; Nucl. Phys. {\bf B431} (1994)
484. }

\lref\ad{P.C. Argyres and M.R. Douglas, ``New phenomena in $SU(3)$
supersymmetric gauge theory," hep-th/9505062; Nucl. Phys. {\bf
B448} (1995) 93.}
\lref\ds{M.R. Douglas and S.H. Shenker, ``Dynamics of $SU(N)$ supersymmetric
gauge theory," hep-th/9503163; Nucl. Phys. {\bf B447} (1995) 271.}

\lref\giveonkut{A. Giveon and D. Kutasov,
``Brane Dynamics and Gauge Theory,'' hep-th/9802067}

\lref\vw{C. Vafa and E. Witten,
``A strong coupling test of $S$-duality,'' hep-th/9408074; Nucl.
Phys. {\bf B431} (1994) 3.}

\lref\monopole{E. Witten, ``Monopoles and
four-manifolds,''  hep-th/9411102; Math. Res. Letters {\bf 1}
(1994) 769.}

\lref\witteni{E. Witten, ``On $S$-duality in abelian
gauge theory,'' hep-th/9505186; Selecta Mathematica {\bf 1} (1995)
383.}

\lref\wittk{E. Witten, ``Supersymmetric Yang-Mills theory
on a four-manifold,''  hep-th/9403193; J. Math. Phys. {\bf 35}
(1994) 5101.}

\lref\lns{A. Losev, N. Nekrasov, and S. Shatashvili, ``Issues in
topological gauge theory," hep-th/9711108; Nucl. Phys. {\bf B 534} (1998)
549. ``Testing Seiberg-Witten
solution," hep-th/9801061.}

\lref\mm{M. Mari\~no and G. Moore, ``Integrating over the Coulomb branch in
${\cal N}=2$ gauge theory," hep-th/9712062;  Nucl. Phys. Proc.
Suppl. {\bf b68} (1998) 336.}

\lref\matone{M. Matone, ``Instantons and recursion relations in ${\cal N}=2$
supersymmetric gauge theory," hep-th/9506102; Phys. Lett. {\bf
B357} (1995) 342.}

\lref\humphreys{J.E. Humphreys, {\it Introduction to Lie algebras and
representation theory}, Springer-Verlag, 1972.}

\lref\fm{R. Friedman and J.W. Morgan,
``Algebraic surfaces and Seiberg-Witten invariants,''
alg-geom/9502026; J. Alg. Geom. {\bf 6} (1997) 445. ``Obstruction
bundles, semiregularity, and Seiberg-Witten invariants'',
alg-geom/9509007.}

\lref\morganbk{J.W. Morgan, {\it The Seiberg-Witten equations and applications
to the topology of smooth four-manifolds}, Princeton University
Press, 1996.}

\lref\DoKro{S.K.~ Donaldson and P.B.~ Kronheimer,
{\it The Geometry of Four-Manifolds}, Clarendon Press, Oxford,
1990.}

\lref\fmbook{R. Friedman and J.W. Morgan,
{\it Smooth four-manifolds and complex surfaces}, Springer Verlag,
1991.}

\lref\tqft{E. Witten,
``Topological Quantum Field Theory,'' Commun. Math. Phys. {\bf 117}
(1988) 353.}

\lref\apsw{P.C. Argyres, M.R. Plesser, N. Seiberg, and E. Witten, ``New ${\cal
N}=2$ superconformal field theories in four dimensions,"
hep-th/9511154; Nucl. Phys. {\bf B461} (1996) 71.}

\lref\eguchisc{T. Eguchi, K. Hori, K. Ito, and S.K. Yang, ``Study of ${\cal
N}=2$ superconformal field theories in four dimensions,"
hep-th/9603002; Nucl. Phys. {\bf B471} (1996) 430.}

\lref\BSV{M. Bershadsky, V. Sadov, and
C. Vafa, ``D-Branes and Topological Field Theories,'' Nucl. Phys.
{\bf B463} (1996) 420; hep-th/9511222.}

\lref\nsc{M. Mari\~no and G. Moore, ``Donaldson invariants for nonsimply
connected manifolds," hep-th/9804104.}

\lref\munnsc{V. Mu\~noz, ``Basic classes for four-manifolds
not of simple type,'' math.DG/9811089.}

\lref\freedman{D. Anselmi, J. Erlich, D.Z. Freedman,
and A.A. Johansen, ``Nonperturbative formulas for central functions
of supersymmetric gauge theories,'' hep-th/9708042;
Nucl.Phys.B526:543-571,1998 ; ``Positivity constraints on anomalies
in supersymmetric gauge theories,'' hep-th/9711035;
Phys.Rev.D57:7570-7588,1998}

\lref \pidtyurin{V. Pidstrigach and A. Tyurin, ``Localization of the Donaldson
invariants along Seiberg-Witten classes,'' dg-ga/9507004.}

\lref\feehan{P.M.N. Feehan and  T.G. Leness, ``$PU(2)$ monopoles and relations
between four-manifold invariants," dg-ga/9709022. }

\lref\baryon{J.M.F. Labastida and M. Mari\~no, ``Twisted baryon
number in ${\cal N}=2$ supersymmetric QCD," hep-th/9702054; Phys. Lett.
{\bf B 400} (1997) 323.}

\lref\seiberg{N. Seiberg, Phys. Lett. {\bf B 318} (1993) 469; Phys. Rev.
{\bf D 49} (1994) 6857.}
\lref\ils{K. Intriligator, R.G. Leigh and N. Seiberg, Phys. Rev. {\bf D 50}
(1994) 1052.}
\lref\is{K. Intriligator and N. Seiberg, ``Lectures on ${\cal N}=1$
supersymmetric gauge theories
and electric-magnetic duality," hep-th/9509066; Nucl. Phys. Proc.
Suppl. {\bf B45} (1996)1. }

\lref\gompfbook{R. E. Gompf and A.I. Stipsicz, {\it An introduction to
four-manifolds and
Kirby calculus}, available in R.E. Gompf's homepage at
http://www.ma.utexas.edu.}

\lref\bilalfer{A. Bilal and F. Ferrari, ``The BPS spectra and
superconformal points in massive
${\cal N}=2$ supersymmetric QCD", hep-th/9706145; Nucl. Phys. {\bf
B516} (1998) 175.}

\lref\ky{H. Kanno and S.-K. Yang, ``Donaldson-Witten function of
massless ${\cal N}=2$ supersymmetric QCD," hep-th/9806015; Nucl. Phys.
{\bf B 535} (1998) 512.}

\lref\kub{T. Kubota and N. Yokoi, ``RG flow near the superconformal
points in ${\cal N}=2$ supersymmetric gauge theories," hep-th/9712054;
Prog. Theor. Phys. {\bf 100} (1998) 423.}

\lref\hr{M. Mari\~no and G. Moore, ``The Donaldson-Witten function for
gauge groups of rank larger than one," hep-th/9802185.}

\lref\bpv{W. Barth, C. Peters and A. Van de Ven, {\it Compact complex
surfaces}, Springer, 1984.}

\lref\beau{A. Beauville, {\it Complex algebraic surfaces}, Cambridge
University Press, 1996.}

\lref\gompfnc{R.E. Gompf, ``A new construction
of symplectic manifolds," Ann. of Math. {\bf 142} (1995) 527.}

\lref\taubes{C.H. Taubes, ``The Seiberg-Witten
invariants and symplectic forms," Math. Res. Letters
{\bf 1} (1994) 809; ``SW $\Rightarrow$ Gr: From the
Seiberg-Witten equations to pseudo-holomorphic curves,"
J. Amer. Math. Soc. {\bf 9} (1996) 845.}

\lref\ko{D. Kotschik, ``The Seiberg-Witten invariants of
symplectic four-manifolds [after C.H. Taubes],"
S\'eminaire Bourbaki, 48\`eme ann\'ee, 1995-6, n. 812. }

\lref\szabo{Z. Szab\'o, ``Simply-connected irreducible four-manifolds with
no symplectic structures," Invent. Math. {\bf 132} (1998) 457. }

\lref\fsrat{R. Fintushel and R.J. Stern, ``Rational blowdowns of smooth
four-manifolds," alg-geom/9505018, J. Diff. Geom. {\bf 46} (1997) 181; ``Knots,
links, and
four-manifolds," dg-ga/9612004; ``Constructions of
smooth four-manifolds," Doc. Math., extra volume ICM 1998, p. 443.}

\lref\brussee{R. Brussee, ``The canonical class and the $C^{\infty}$
properties of K\"ahler surfaces,'' alg-geom/9503004. }

\lref\kmone{P.B. Kronheimer and T.R. Mrowka, ``Embedded
surfaces and the structure of Donaldson polynomials," J. Diff. Geom. {\bf
41} (1995) 573. }

\lref\kmad{P.B. Kronheimer and T.R. Mrowka,
``The genus of embedded surfaces in the projective plane," Math. Res. Lett.
{\bf 1} (1994) 797.}

\lref\jpark{J. Park, ``The geography of irreducible four-manifolds,"
Proc. Amer. Math. Soc. {\bf 126} (1998) 2493.}

\lref\mms{J.W. Morgan, T.R. Mrowka, and Z. Szab\'o, ``Product formulas along
$T^3$ for Seiberg-Witten invariants," Math. Res. Lett. {\bf 4} (1997) 915.}

\lref\mszabo{J.W. Morgan and Z. Szab\'o, ``Embedded
genus-$2$ surfaces in four-manifolds,"
Duke Math. Journal {\bf 89} (1997) 577.}

\lref\fsonly{R. Fintushel and R.J. Stern, ``Nonsymplectic
four-manifolds with one basic class,"  preprint.}

\lref\hh{F. Hirzebruch and H. Hopf, ``Felder von
Fl\"achenelementen in $4$-dimensionalen Mannigfaltigkeiten",  Math. Ann.
{\bf 136} (1958) 156.}

\lref\lloz{J.M.F. Labastida and C. Lozano, ``Duality
symmetry in twisted ${\cal N}=4$ supersymmetric gauge theories in four
dimensions," hep-th/9806032.}

\lref\gr{M. Mari\~no and G. Moore, ``The Donaldson-Witten
function for gauge groups of rank larger than one," hep-th/9802185.}

\lref\na{J.M.F. Labastida and M. Mari\~no, ``Non-abelian monopoles on
four-manifolds,"
hep-th/9504010; Nucl. Phys. {\bf B 448} (1995) 373. ``Polynomial invariants for
$SU(2)$ monopoles,"
hep-th/9507140; Nucl. Phys. {\bf B 456} (1995) 633.}

\lref\tqcd{S. Hyun, J. Park and J.S. Park, ``Topological QCD,"
hep-th/9503201; Nucl. Phys. {\bf B 453} (1995) 199.}

\lref\afre{D. Anselmi and P. Fr\'e, ``Topological sigma
model in four-dimensions and triholomorphic maps," hep-th/9306080; Nucl.
Phys. {\bf B 416}
(1994) 255. ``Gauge hyperinstantons and monopole equations,"
hep-th/9411205; Phys. Lett. {\bf B 347}
(1995) 247.}

\lref\tmfour{M. Alvarez and J.M.F. Labastida, ``Topological
matter in four dimensions," hep-th/9404115; Nucl. Phys. {\bf B 437} (1995)
 356.}

\lref\oko{C. Okonek and A. Teleman, ``Recent developments in
Seiberg-Witten theory and complex geometry," alg-geom/9612015, and
references therein.}

\lref\pt{V. Pidstrigach and A. Tyurin, ``Localisation of the Donaldson
invariants along Seiberg-Witten basic classes," dg-ga/9507004.}

\lref\fele{P.M.N. Feehan and  T.G. Leness, ``$PU(2)$ monopoles and relations
between four-manifold invariants," dg-ga/9709022; ``$PU(2)$
monopoles I: Regularity, Uhlenbeck compactness, and
transversality," dg-ga/9710032; ``$PU(2)$ monopoles II:
Highest-level singularities and relations between four-manifold
invariants," dg-ga/9712005.}

\lref\eqc{J.M.F. Labastida and M. Mari\~no, ``Twisted ${\cal N}=2$
supersymmetry
with central charge and equivariant cohomology," hep-th/9603169; Comm.
Math. Phys.
{\bf 185} (1997) 323.}

\lref\GrHa{P.~ Griffiths and J.~ Harris, {\it Principles of
Algebraic Geometry},  J.Wiley and Sons, 1978. }

\lref\mmth{M. Mari\~no, ``The geometry of supersymmetric
gauge theories in four dimensions," hep-th/9701128.}

\lref\mmp{M. Mari\~no, G. Moore and G. Peradze, ``Superconformal
invariance and the geography of four-manifolds," hep-th/9812055.}

%%%%%%%%%%

%%%%%%%%%%
\Title{\vbox{\baselineskip12pt
\hbox{YCTP-P30-98 }
\hbox{math.DG/9812042}
}} {\vbox{\centerline{   Four-Manifold Geography   }
\centerline{}
\centerline{ and }
\centerline{}
\centerline{ Superconformal Symmetry }}
}
\centerline{Marcos Mari\~no, Gregory Moore, and Grigor Peradze}

\bigskip
{\vbox{\centerline{\sl Department of Physics, Yale University}
\vskip2pt
\centerline{\sl New Haven, CT 06520, USA}}

\centerline{ \it marcos.marino@yale.edu }
\centerline{ \it moore@castalia.physics.yale.edu }
\centerline{ \it grigor.peradze@yale.edu }

\bigskip
\bigskip
\noindent
A compact  oriented 4-manifold is
defined  to be of  ``superconformal simple type''
if certain   polynomials in the basic
classes (constructed using the Seiberg-Witten
invariants) vanish identically.
We show that all known 4-manifolds
of $b_2^+>1$ are of superconformal
simple type, and that the numerical
invariants of 4-manifolds of
superconformal simple type satisfy
a generalization of the Noether inequality.
We sketch how
these phenomena  are predicted by the
existence of certain four-dimensional
superconformal quantum field theories.

\Date{Dec. 6, 1998}

%\draft

\newsec{Introduction}

On several occasions  insights
into the physics of quantum field
theory and string theory  have suggested new
results and techniques in mathematics, particularly
in   geometry and topology.
One of the main focal points of activity
in physics in the past few years
has been the study of superconformal
field theories in four dimensions. In this letter
we show that some of the recent progress in
physics leads to new results on the geography
of four-manifolds.

This letter is primarily intended for a mathematical
audience. The physical motivation for  our
results is summarized briefly in section five.
A more extensive account of this work
can be found in \mmp.

\newsec{ SST Manifolds}

Let $X$ be a compact, oriented four-manifold
with $b_2^+>1$. In this letter we will address
the relation between the classical numerical
invariants of $X$:
\eqn\vari{
\eqalign{
\chi_h & := { \chi + \sigma \over  4}
= { 1- b_1 + b_2^+ \over  2}, \cr
c_1^2  & := 2 \chi + 3 \sigma = 4 - 4b_1 + 5 b_2^+-b_2^- ,  \cr}}
and its Seiberg-Witten invariant, which
is a map $SW: {\rm Spin^c}(X) \rightarrow \IZ$ from
the ${\rm Spin^c}$ structures of $X$ to the integers.
$SW$ is a topological invariant  of $X$
defined using a signed sum over solutions
to the monopole equations \monopole.
We will identify ${\rm Spin^c}(X) $ with
elements of $H^2 (X, \IZ)$ congruent to $w_2 (X)$ mod $2$,
and say that   $x\in {\rm Spin^c}(X) $ is a
basic class if $SW(x) \not=0$.
In this letter we
 will assume that $X$ is of (Seiberg-Witten) simple type, {\it
i.e.},   if $x^2 \not=c_1^2$ then
$SW(x)=0$. If $X$ is of
simple type and the SW invariant is not trivial, then $X$ must
obey the Noether condition, {\it i.e.} $\chi_h\in \IZ$.

Choose an integral lifting $\upsilon$ of $w_2(X)$ and
consider  the
twisted SW series:
\eqn\swseries{
SW_X^{w_2 (X)} (z) := \sum_{x} (-1)^{\upsilon^2 + \upsilon\cdot x \over 2} SW
(x) {\rm e}^{z  x}.}
This is a finite sum \monopole.  A change of lifting $\upsilon\rightarrow
\tilde \upsilon$
alters $SW_X^{w_2 (X)}$  by a sign
$(-1)^{\bigl({\upsilon-\tilde \upsilon \over 2}\bigr)^2}$.

\bigskip
\noindent
{\bf Definition 2.1. } Let $X$ be a compact, oriented 4-manifold of simple type
with
$b_2^+>1$. We say that ``$X$   is  SST'' if
$SW_X^{w_2 (X)} (z)$ has a zero at $z=0$ of order $\ge \chi_h - c_1^2 -3$.

The phrase ``$X$   is  SST'' is short for
``$X$ is of   superconformal simple
type.'' The terminology comes from the physical motivation discussed at
the end of this letter.
 The SST condition is equivalent to the condition
that either $c_1^2 - \chi_h + 3 \geq 0$ or
\eqn\sumrules{
\sum_{x}(-1)^{\upsilon^2 + \upsilon \cdot x \over 2} SW (x)
x^{k} =0, \,\,\,\,\,\,\,\ k=0, \dots,  \chi_h - c_1^2 -4, }
where $x^k$ is naturally understood as an element in
$\big( {\rm Sym}^k (H^2(X, \IZ)\bigr)^*$ acting through the
intersection form. Notice that, if $\chi_h + \sigma$ is even (odd),
the expressions of the
form \sumrules\
with $k$ odd (even) are automatically zero. This is easily proved
using the fact that, if $x$ is a basic class, then so is $-x$, and
\eqn\wittfor{
(-1)^{\upsilon^2 - \upsilon\cdot x \over 2}SW(-x) = (-1)^{\chi_h + \sigma}
 (-1)^{\upsilon^2 + \upsilon\cdot x \over 2} SW (x),}
as one easily checks using the Wu formula and
the   behavior of the SW invariants under the
involution $x \rightarrow -x$ \monopole. Therefore,
there are only  $[\half(\chi_h- c_1^2)] -1$ nontrivial
equations in \sumrules, where
$[ \cdot ] $ is the greatest integer function.

\newsec{Two theorems about SST manifolds}

In this section, we will use the known behavior of the
Seiberg-Witten invariants to
state some properties of SST  manifolds. The detailed
proofs of Theorems 3.1 and 3.2  can be found
in \mmp. We only make some
brief remarks on the proofs here.

\bigskip
\noindent
\thm{3.1} The property ``$X$ is   SST''  is preserved by
blowup, fiber sum along c-embedded tori,   knot surgery,
and generalized log transforms.
\bigskip
The proof of this theorem is a simple consequence of
the  known behavior of
the SW invariants under these operations \fsrat\mms. The blowup is particularly
instructive. If $\widehat X
=X \sharp \overline{ \IC P^2}$ denotes the blownup manifold,
with exceptional divisor $E$, one has
\eqn\swsbup{
SW_{\widehat X} ^{w_2 (\widehat X)} (z) = -2 \sinh (z  E) SW_X^{w_2
(X)} (z),}
where we choose the lift of $w_2 (\widehat X)$ to be $\hat
\upsilon=\upsilon+E$.
Using $c_1^2 (\widehat X) = c_1^2 (X)-1$, it   immediately
follows that if $X$ is SST then $\widehat X$ is SST. This computation shows
that the sign factor in \swseries\ is crucial.

The proof for the fiber sum along c-embedded tori is based on the
gluing formulae of \fsrat\mms. Using these formulae, one shows that,
if $X=X_1 \sharp_{T_1=T_2} X_2$ is the fiber sum of two SST manifolds,
then
\eqn\twistser{
SW_X^{w_2 (X)}(z) = \pm 4 \,(\sinh zT)^2 \,SW_{X_1}^{w_2 (X_1)}
(z) \cdot SW_{X_2}^{w_2 (X_2)} (z),}
where the $\pm$ sign depends on the choice of the liftings
and $T$ is the common torus. The proof that
$X$ is SST is easy arithmetic based on \twistser.

The case of knot surgery also  follows from \fsrat.
Finally, we consider the
  generalized log transform for manifolds with a cusp neighborhood. This is
done by   doing $p-1$ blowups,
$X \rightarrow X \sharp
(p-1){\overline{ \IC P^2}}$, and then performing a rational blowdown along
a $C_p$ configuration. See \fsrat\jpark\gompfbook\ for details. Again, the
proof is simple given   the known behavior of the
SW invariants under these operations.

\thm{3.2} If $X$ is a compact complex surface of
$b_2^+>1$ then $X$ is SST.
\bigskip

\noindent
{\it Proof}:
Since property SST is preserved by blowup,
we can restrict to minimal complex surfaces. Using the
Kodaira-Enriques classification we need only
consider   elliptic surfaces and
surfaces of general type. For minimal surfaces of general type,
the property
SST is a consequence of the Noether inequality. For
minimal elliptic fibrations with no multiple fibers
 the computation of the twisted series
\swseries\ follows from  \fm\brussee\ and is:
\eqn\holef{
SW^{w_2(X)}_X (z) =   \bigl(2 \sinh (z f)) ^{\chi_h
+ 2g-2},}
Here $g\geq 0$ is the genus of the base of the fibration
and $f$ is the class of the fiber.
Taking into account that $c_1^2=0$ the SST property
follows immediately. In this computation we used
the lift $\upsilon=c_1(K)$, with $K$ the canonical bundle of the
elliptic surface. Since the SST property is preserved by generalized
log transforms  it follows that any minimal elliptic
surface is SST. $\spadesuit$

Most  of the recent constructions of ``exotic" manifolds
(symplectic but noncomplex
manifolds or irreducible nonsymplectic manifolds)
are based on the constructions we have
considered, and use complex surfaces as their building blocks. Therefore, they
are all SST. For example, the
  SST property is also preserved by rational blowdowns of tautly
embedded $C_{n-2}$ configurations in the simply connected elliptic fibrations
$E(n)$ (for $n\ge 4$). In particular, after a rational blowdown along one
$C_{n-2}$
in $E(n)$, one obtains an interesting exotic
manifold $Y(n)$. Using  \fsrat\ one
 finds that  $c_1^2=\chi_h-3$ for
$Y(n)$, so it is  SST.
Based on the above results and the examination of many other examples, we
accordingly formulate the  following

\bigskip
\noindent
{\bf Conjecture 3.3}. Every compact, oriented 4-manifold of
$b_2^+(X)>1$  is SST.
\bigskip

We will briefly explain the physical reasons for this conjecture in
section 5. It would be desirable to have more results
along the lines of this section for other constructions of
four-manifolds, such as arbitrary rational blowdowns or
fiber sums along surfaces of genus $g>1$. This may help
to provide a mathematical proof of the above conjecture.

\newsec{Bounds on the number of basic classes for
SST manifolds}

One of the most important properties of SST manifolds
is that, if they support any basic classes at all, then
there is a lower bound on the number of such classes
in terms of the classical numerical invariants of $X$.
A corollary of this lower bound is a generalized Noether
inequality that relates the values of the numerical invariants $c_1^2$,
$\chi_h$, to the number of basic classes. This gives a deep relation
between the Seiberg-Witten invariants and the geography of four-manifolds.

Let $\CB_X$ denote the set of basic classes of $X$. If $x$ is a basic class,
then so is $-x$, and it is convenient to mod out by this involution
when counting the number of basic classes. We will then
say that $X$ has $B$ basic classes if the set $\CB_X/\{ \pm 1 \}$ consists of
$B$ elements.

\bigskip
\noindent
{\bf Theorem 4.1}(Generalized Noether inequality). Let $X$ be   SST. If $X$ has
$B$ distinct basic
classes and $B>0$, then
\eqn\bound{
B \ge \biggl[ {\chi_h -c_1^2 \over 2} \biggr].}
In particular, $c_1^2 \ge \chi_h -2B-1$.
\bigskip

The theorem is a direct consequence
of the conditions \sumrules\ for SST manifolds.
The detailed proof
can be found in \mmp.
The bound \bound\ is in fact
sharp since it is saturated by the simply-connected elliptic
fibrations $E(n)$.
It follows from Theorem 4.1 that
 manifolds with only one basic class
must satisfy $c_1^2 \ge \chi_h -3$. This inequality is saturated, as we have
noted, by the exotic
manifolds $Y(n)$ obtained by rational blowdown of $E(n)$ \fsrat.
All known constructions, e.g. those found in \fsrat,
 of manifolds with one basic class satisfy
 our inequality.

Given a metric on
$X$ one can also give an upper bound on $B$ in
terms of certain Riemannian functionals,
using the ideas in \monopole.  Thus, there exist
topological lower bounds on these Riemannian
functionals.

\newsec{Sketch of the relation to superconformal field theory}

The above results and definitions are motivated
by the existence of superconformal fixed points
in certain physical theories. In this section, we summarize
the main physical ideas behind the above results.
Further discussion can be found in
\mmp.

Donaldson theory can be interpreted as a ``topologically
twisted'' version
of an $\CN=2$ supersymmetric Yang-Mills theory \tqft. The twisting
procedure, which can be understood as a redefinition of
the coupling to gravity, produces quantum field theories
which are formally ``topological," in the sense that the correlation
functions do not depend on the metric of the four-manifolds. In this way,
certain correlation functions of the
twisted ${\cal N}=2$ theory
are identified with Donaldson invariants. These correlation
functions can be assembled into a generating function
known as the Donaldson-Witten partition function
and denoted by  $Z_{DW}$.

The data specifying the
  general $d=4,\CN=2$ supersymmetric gauge theory is
a choice
of a compact Lie group $G$ (the gauge group), a finite dimensional
representation $R$ (the matter hypermultiplet representation) and a
$G$-invariant tensor
in $R \otimes R^*$ (the masses of the matter hypermultiplets).
The
procedure of ``topological twisting''
only depends on the structure of the ${\cal N}=2$
supersymmetry algebra, so any of these models can be twisted to
obtain a topological field theory. A particularly interesting
case is the topological theory defined by the gauge group $SU(2)$,
with one matter multiplet in the fundamental
representation of $SU(2)$, denoted ${\bf 2}$. This theory
defines the moduli space associated to the nonabelian monopole equations,
and it can be analyzed following the lines of Donaldson theory.
(See \mmth\ for a review from the point of
view of topological field theory, and \feehan\ for a recent review
of rigorous results.) In particular, one can formally define
topological invariants using intersection theory in the moduli
space of solutions. An interesting aspect of this
theory is that the mass term can be interpreted as the formal
parameter of the equivariant cohomology of  the moduli space  with respect
to a $U(1)$ group action \eqc.
Therefore, the correlation functions
of this topological theory, considered as   functions of the mass, formally
compute the  equivariant
intersection theory on the moduli space of nonabelian monopoles. More
generally, for the theory specified
by $G,R$ the masses $m$ are the parameters for the
equivariant cohomology of a group action on the
moduli space of the corresponding (nonabelian)
monopole equations.

%In particular, in the example of  theories
%``with $N_f\ge 1$ matter
%multiplets,'' i.e., with
%$R= N_f {\bf 2} $  the global symmetry is $U(1)^{N_f}$ and the
%parameters of the corresponding equivariant cohomology are %a collection of
% masses
%$m_\alpha$, $\alpha=1, \cdots, N_f$.
%

The low-energy dynamics of the $SU(2)$ ${\cal N}=2$ supersymmetric gauge
theories
was solved by Seiberg and Witten in \swi\ in
terms of certain families of elliptic curves.
 In the case when there are
no matter multiplets  Witten
used this to obtain  the topological
correlation functions, {\it i.e.}  the Donaldson
invariants,
in terms of the SW invariants \monopole. The result
of \monopole\ can be
generalized to the
theories with matter, as explained in \mw, and one obtains a generalization of
Witten's formula, {\it i.e.}, an exact expression for the
generating function of equivariant
intersection
numbers,   $Z_{DW}(p,S;m)$.
In this function, $p,S$ correspond to (equivariant) differential forms of
degree
four and two, respectively, on the corresponding moduli space, and generalize
the usual observables of Donaldson-Witten theory.

We now describe the function $Z_{DW}(p,S;m)$ for the
case $G=SU(2), R= {\bf 2}$ in more detail.
The Seiberg-Witten family of elliptic curves  is
parametrized by $(u,m)\in \IC^2$ and given by:
\eqn\swc{
y^2 = x^2(x-u) + 2m x -1 .
}
The curve is easily put into standard form
$y^2 = 4x^3 - g_2 x - g_3$, with
$g_2(u;m) = {4 \over  3} (u^2-6m)$,
$g_3(u;m) = {1 \over  27}(8 u^3 -72 m u + 108)$, and
discriminant $\Delta(u;m) = g_2^3 - 27 g_3^2$.
This discriminant is a cubic in $u$ and has
three roots $u_j(m)$, $j=1,2,3$. For generic,
but fixed, values of $m$
one of the periods of \swc\ goes to infinity
as $u \rightarrow u_j$ while the other period,
$\varpi_j\equiv \varpi(u_j(m);m)$ remains finite, and in fact is given
by $(\varpi_j )^2 = g_2/(36 g_3)$. The Donaldson-
Witten partition function on a manifold $X$
of simple type and $b_2^+(X)>1$
is   given by a sum over the
singular fibers of the Weierstrass family
\swc\ and over the basic classes of $X$:
\eqn\dwpf{
\eqalign{
Z_{DW}(p,S;m) & = k \sum_{j=1}^3 \biggl({g_2^3(u_j(m);m) \over
\Delta'(u_j(m);m) } \biggr)^{\chi_h} (\varpi_j(m))^{7 \chi_h - c_1^2}\cr
\sum_x  & SW(x) (-1)^{(\upsilon^2 + \upsilon \cdot x)/2}
\exp\biggl[ 2 p u_j + S^2 T_j - i {  (S,x) \over  2\varpi_j  }  \biggr] \cr}
}
Here $\Delta' = { \p \over  \p u} \Delta$,
$T_j=-{1 \over 24} \bigl(
(\varpi_j)^{-2}  -8 u_j \bigr) $,
and $k$ is a nonvanishing constant, independent of
$p,S,m$. Similar expressions hold for gauge groups
$G$ of rank bigger than one \hr.

At first sight, \dwpf\ and its generalizations are
disappointing, since they imply that the more general topological
theories and moduli equations obtained from general supersymmetric gauge
theories do not contain new topological information.
While \dwpf\ and its generalizations seem
disappointing, defeat can
be turned into victory by using known properties
of the  physical theories to learn about the SW invariants.
{}From the mathematical viewpoint, the parameters $m$
of equivariant cohomology are formal expansion
parameters, but from the physical viewpoint
it is clear that we should regard
$Z_{DW}(p,S;m)$
as nontrivial locally analytic functions of   $m$.
At certain points of the moduli space, i.e., at
special values of $m$, the low-energy SW
theory is in fact a nontrivial
superconformal
field theory ({\it i.e.} a conformal field theory with ${\cal N}=2$
supersymmetry) \apsw. We now describe how
the  analytic  structure  of $Z_{DW}(p,S;m )$
at such points can lead to the nontrivial predictions
for topology described above.

Let us summarize some facts about the
SW curve for $G=SU(2)$, $R={\bf 2}$.  Substituting
$m= {3 \over  2} + z, u= 3 + 2z + \delta u$
into \swc\ and taking $z, \delta u$ small leads to the
curve
\eqn\cusp{
y^2 = x^3 - 2z x - \delta u .
}
The family \cusp\
  develops a cusp singularity as $z,
\delta u \rightarrow 0$.
Indeed, when  $z \rightarrow 0$, two of the roots of
$\Delta(u;m)=0$, call them $u_\pm(m)$, coincide.
By scaling one finds
$\delta u_{\pm} \sim z^{3/2}$. In the limit $z \rightarrow 0$
the period $\varpi_\pm $ diverges as $z^{-1/4}$, while
$g_2(u_\pm(m);m) \sim z$  and $\Delta'(u_\pm(m);m)
\sim \delta u_{\pm} \sim z^{3/2}$. At the third
singularity all the various factors in \dwpf\ are
given by nonvanishing analytic series in $z$, but,
evidentally, the contributions from $u_\pm(m)$ contain
factors which are diverging or vanishing as $z \rightarrow 0$.
What can we say about the behavior of the
complete function $Z_{DW}$ as $z \rightarrow 0$?

The detailed expansion of $Z_{DW}$ as a power series
in $z$ is discussed at length in \mmp. An easy
consequence of \wittfor\ is that $Z_{DW}(z)$ is
a Laurent series in integral powers of $z$.
Using the facts from the previous paragraph one
checks that the Laurent series \dwpf\ has at most a
finite order pole   of order  $z^{(c_1^2 - \chi_h)/4}$.
Thus, if $c_1^2- \chi_h<0$ there could
be a pole in $z$.  However, the
physical origin of $Z_{DW}$ leads to a
powerful and well-founded principle:

{\it If $X$ is compact, and the moduli space of solutions
to the relevant abelian
monopole equations is compact, then
$Z_{DW}$ must be a regular analytic function of $z$
near $z=0$.}

The reason for this principle is that divergences in
topological field theory, even at superconformal
points, can only arise from infrared divergences in
spacetime or in moduli space. Since $X$ is compact
there are no divergences in spacetime. If the
moduli space is compact then $Z_{DW}$ cannot diverge.

If $c_1^2 - \chi_h < -3$
then  the coefficients of the would-be poles in
the $z$-expansion
are polynomials in $p, S^2$, and $(S,x)$. The regularity
of $Z_{DW} $ at $z=0$ thus implies nontrivial
constraints on these quantities. Unfortunately,
the detailed expansion described in \mmp\
is somewhat complicated.
Nevertheless, since $\varpi_\pm \sim z^{-1/4}$
it follows from expanding the exponential in
\dwpf\ that if the sum rules \sumrules\ are
satisfied then $Z_{DW}$ is indeed analytic at
$z=0$.
This sufficient condition is the SST condition.
As we have checked above, all known 4-manifolds
of $b_2^+>1$ satisfy the condition.

Close examination of the Laurent expansion
in $z$ of \dwpf\ reveals that  the SST property
is sufficient, but not necessary,
for regularity  \mmp. It is thus logically
possible that there exist 4-manifolds of
simple type which are not SST, and yet
do not violate our physical principle of regularity.
However, the full statement of conditions for
regularity of $Z_{DW}$ is rather involved,
and we find this possibility extremely unlikely, both
from the point of view of quantum field theory
and from the point of view of four-manifold topology.
It is worth noting that, even if we do drop the
assumption that $X$ is SST, regularity of
$Z_{DW}$ implies that if
$X$ has one basic class then  it
satisfies the generalized Noether inequality
$c_1^2 \geq \chi_h -3$ \mmp.

One can probably even relax the condition that
$X$ is of simple type. Every 4-manifold of
$b_2^+>1$ is of generalized simple type \mw.
Using equations $(7.12)$ and $(11.27)$ of
\mw\ it should be possible to extend the
above results to those (hypothetical) manifolds
of $b_2^+>1$ which are not of simple type.

It would be very interesting to see whether the
reasoning we have described here leads to further
results in topology using other
$\CN=2$ superconformal fixed points.

\bigskip
\centerline{\bf Acknowledgements}\nobreak
\bigskip

We would like to thank E. Calabi,
D. Freed, R. E. Gompf, T.J. Li and E. Witten for very
helpful discussions and correspondence.
This work  is supported by DOE grant
DE-FG02-92ER40704.

\listrefs

\bye